\documentclass[a4paper,10pt,leqno,dvipdfmx]{amsart}
\usepackage{amsmath, amssymb, amscd, amsthm}
\usepackage{mathrsfs}
\usepackage{graphicx}
\usepackage{enumerate}

\makeatletter

\def\@settitle{\begin{center}
  \baselineskip14\p@\relax
  \normalfont\LARGE\bfseries
  \@title
  \ifx\@subtitle\@empty\else
     \\[1ex] 
     
     \normalsize\mdseries\@subtitle
  \fi
 \ifx\@didication\@empty\else
     \\[2ex] 
     
     \large\mdseries\it\@dedication
  \fi
  \end{center}
}
\def\subtitle#1{\gdef\@subtitle{#1}}
\def\@subtitle{}
\def\dedication#1{\gdef\@dedication{#1}}
\def\@dedication{}
\makeatother

\makeatletter
\renewcommand{\section}{\@startsection
{section}{1}{0mm}{5mm}{2mm}{\raggedright\bfseries}}

 \@addtoreset{equation}{section}
\makeatother

\title{Leaf as a Poincar\'e convex domain associated with an endomorphism on a real inner product space}

\author{Hiroyuki Ogawa}

\address{
Department of Mathematics, 
Graduate School of Science, 
Osaka University, 
Toyonaka, Osaka 560-0043, Japan}
\email{ogawa@math.sci.osaka-u.ac.jp}

\theoremstyle{plain}

\theoremstyle{definition}

\begin{document}

\newcommand\vc[1]{{\mathsf{#1}}}
\newcommand\Proj{{\mathbb P}}
\newcommand\Hu{{\mathfrak{H}}}
\newcommand\D{{\mathcal D}}
\newcommand\N{{\mathbb N}}
\newcommand\C{{\mathbb C}}
\newcommand\Z{{\mathbb Z}}
\newcommand\R{{\mathbb R}}
\newcommand\Q{{\mathbb Q}}
\newcommand\M{\mathrm{M}}
\newcommand\End{\mathrm{End}}
\newcommand\EndR{\End}
\newcommand\id{{\text{id}}}
\newcommand\pvs[1]{{{\mathbb L}_{{#1}}}}
\newcommand\pv{\pvs{\varphi}}
\newcommand\pvk{\pvs{\varphi+k\id_V}}
\newcommand\pvj[1]{\pvs{J(0:{#1})}}
\newcommand\real{{\mathrm{Re}\,}}
\newcommand\imag{{\mathrm{Im}\,}}
\newcommand\ssqrt[1]{{\sqrt{{#1}}}}
\newcommand\unitCpx{\ssqrt{-1}}
\newcommand\mt[1]{{{\mathsf #1}}}
\newcommand\mtE{{\mt E}}
\newcommand\mtJ{{\mt J}}
\newcommand\mtK{{\mt K}}
\newcommand\mtS{{\mt S}}
\newcommand\gC{{C_{\rm P}}}

\newcommand\tp{{{}^t}}
\newcommand\trans{\tp}
\newcommand\tpp{{\tp\!}}

\newcommand\lmat{\left(}
\newcommand\rmat{\right)}

\newcommand\matrixformnull[2]{
  \begin{array}{@{\,}#1@{\,}}#2\end{array}}
\newcommand\matrixform[2]{\lmat\matrixformnull{#1}{#2}\rmat}

\newcommand{\thmenv}[2]{
  \thmnumbegin
  \begin{#1}
    #2
  \end{#1}
  \thmnumend
}

\newcounter{ThmNum}
\newcommand{\thmnumbegin}{\setcounter{ThmNum}{0}}
\newcommand{\thmnumend}{}
\renewcommand{\theThmNum}{\alph{ThmNum}}
\newcommand{\thmNum}{\ \hspace{-1em} \addtocounter{ThmNum}{1}\hbox{\ \ {\rm({\theThmNum})}} \, }
\newcommand{\setThmNum}[1]{\setcounter{ThmNum}{#1}}

\newcommand\subFigure[3]{\begin{minipage}[t]{#1}
	\small
	\centering
	#3 \\
	#2

	\end{minipage}
}

\begin{abstract}
We define a subset of the closure of the upper half plane associated with an endomorphism on a real inner product space, which is called the leaf. 
When the dimension of the space is at least 3, the leaf is a convex with respect to the Poincar\'e metric, and contains all eigenvalues with nonnegative imaginary part. 
Moreover, the leaf of a normal endomorphism is the minimum Poincar\'e convex domain containing all eigenvalues with nonnegative imaginary part. 
The most commonly studied convex domain containing eigenvalues is number range. 
Numerical range is convex with respect to the Euclidean metric on $\C$, so numerical range has less information than leaf about real eigenvalues. 
We provide a new visual approach to endomorphisms.
\end{abstract}

\maketitle

\markboth{Hiroyuki Ogawa}
{Leaf of endomorphism}

\section{Introduction}\label{sec:intro}

In this paper we define a subset of the closure of the upper half plane associated with an endomorphism on a real inner product space, which is called the leaf, and show the structure theorems (Theorem \ref{convex}, \ref{leaf-eigen} and \ref{normal-polygon}). 
The structure theorems on leaf correspond to the fundamental theorems on numerical range (Theorem \ref{Toe-Hau}, \ref{Win} and \ref{polygon-NR}), which is the well-known convex domain containing all eigenvalues, as below. 
Both numerical range and leaf are simply connected domains containing eigenvalues, but numerical range is convex with respect to the Euclidean metric on $\C$, and leaf is convex with respect to the Poincar\'e metric
on the closure of the upper half plane. 
The structure theorems on leaf and the fundamental theorems on numerical range correspond to each other, but they are proved in the different ways because of the difference in the shape of leaf and numerical range. 
The most significant difference is in the points shared with the real axis, whereas the intersection of numerical range and the real axis is an interval containing 
all real eigenvalues, whereas that of leaf is the set of all real eigenvalues. 
Since leaf hold information on individual real eigenvalues, it is possible to study the behaviour of an endomorphism around the eigenspace with respect to a real eigenvalue. 

Let $V$ be a finite dimensional inner product space over $\R$. 
We denote the ring of endomorphisms ($\R$-linear maps) on $V$ by $\EndR(V)$. 
Since the dimension of the space V is important in this paper, an endomorphism on a $d$-dimensional space V is sometimes simply referred to as an endomorphism of dimension $d$. 
The angle $\theta=\theta(\vc v_1,\,\vc v_2)\in[0,\,\pi]\subset\R$ of two nonzero vectors $\vc v_1$, $\vc v_2\in V^\times$ ($=V\smallsetminus\{\vc 0\}$) is defined by $\vc v_1\cdot\vc v_2=\|\vc v_1\|\,\|\vc v_2\|\,\cos\theta$. 
Let $\varphi$ be an endomorphism on $V$. 
We consider the continuous map $\pv:V^\times\to\C$ is defined by 
\[\pv(\vc v)=\frac{\|\varphi(\vc v)\|}{\|\vc v\|}\,e^{\sqrt{-1}\,\theta(\vc v,\,\varphi(\vc v))} \qquad \text{for $\vc v\in V^\times$.}\] 
When $\varphi(\vc v)=\vc 0$, we don't define the angle $\theta(\vc v,\,\varphi(\vc v))$, but can put $\pv(\vc v)=0$ because of $\|\varphi(\vc v)\|=0$. 
Since $0\le\theta(\vc v,\,\varphi(\vc v))\le\pi$, the imaginary part of $\pv(\vc v)$ is nonnegative. 
$\pv(\vc v)$ belongs to the closure of the upper half plane $\overline{\Hu}=\Hu\cup\R\cup\{\infty\}$, where $\Hu=\{z\in\C\mid \imag z>0\}$. 
\thmenv{definition}
{
The image of $\pv$ 
\[\Psi(\varphi)=\{\pv(\vc v)\mid\vc v\in V^\times\}\subset\overline{\Hu}\]
is called the {\it leaf} of the endomorphism $\varphi$. 
}

\thmenv{example}
{\label{eg-leaf}
Figure 1 contains leaves of five endomorphisms on an real Euclid spaces $\R^n$ with standard inner product. 
Such endomorphism are given by some square matrices and direct sums of square matrices. 
Put $\mtE=\matrixform{rr}{1&0\\0&1}$, 
$\mtJ=\matrixform{rr}{0&1\\{-}1&0}$, 
$\mtK=\matrixform{rr}{0&1\\1&0}$, and 
$J(\lambda;\,d)$ be the Jordan cell with eigenvalue $\lambda\in\R$ and of degree $d$. 
The straight line at the bottom on each figure is the real axis.
\begin{figure}[htb]
 \subFigure{7.5em}{(a) $(1)\oplus(3\,\mtJ{+}2\,\mtK)$ \\ \quad $\oplus({-}4\,\mtE{+}10\,\mtJ)$}{\includegraphics[scale=0.30]{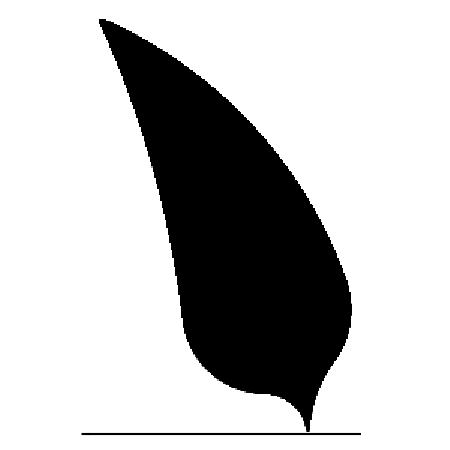}}
 \subFigure{7.5em}{(b) $\matrixform{@{}r@{\ \ }r@{\ \ }r@{}}{0&0&2\\1&0&0\\0&1&0}$}{\includegraphics[scale=0.30]{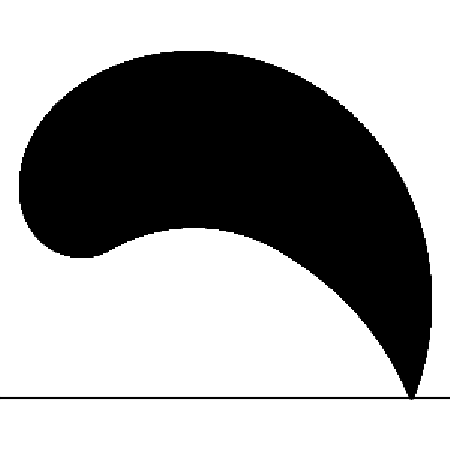}}
 \subFigure{7.5em}{(c) $J(0;2)$}{\includegraphics[scale=0.32]{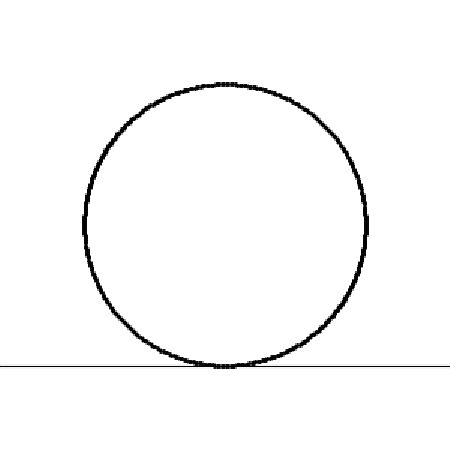}}
 \subFigure{7.5em}{(d) $J(0;3)$}{\includegraphics[scale=0.32]{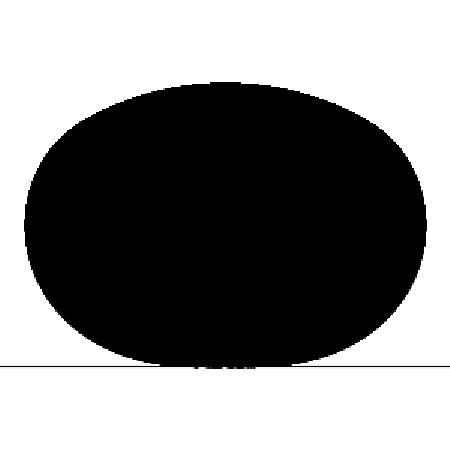}}
 \subFigure{7.5em}{(e) $(-2)\oplus(6\,\mtJ)$ \\ \quad $\oplus(2\,\mtE+\mtJ)$}{\includegraphics[scale=0.30]{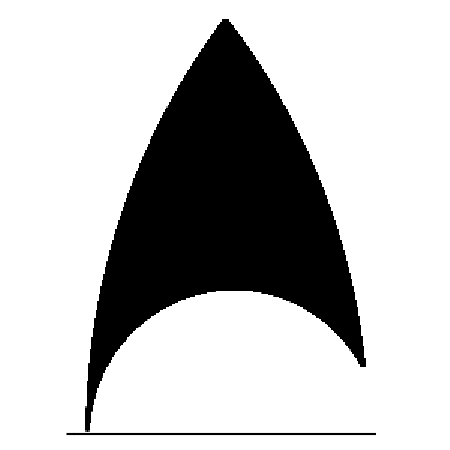}}
 \caption{Examples of leaves}
\end{figure}
}

Leaf is a bounded closed set in the Euclidean metric on $\C$. 
Numerical range is convex with respect to the Euclidean metric on $\C$, but, in generally, leaf is not convex on the Euclidean metric on $\C$. 
In \S 4, we will show that leaf is path-connected (see Proposition 4.2), and that leaf and the real axis share only all real eigenvalues (see Theorem 4.3). 
In \S 6, we will see that leaf of endomorphism of dimension at most $2$ is simple (see Theorem 6.1). 
The purpose of this paper is to show the following structure theorems of leaf of endomorphism of dimension at least $3$. 

\thmenv{mainthm}
{
[Poincar\'e Convexity]
\label{convex}
Leaf is convex with respect to the Poincar\'e metric. 
}

\thmenv{mainthm}
{
[Eigenvalue inclusion]
\label{leaf-eigen}
Leaf contains all eigenvalues whose imaginary part is nonnegative. 
}

The smallest convex domain containing all eigenvalues with nonnegative imaginary part is the convex filled geodesic polygon whose vertices are eigenvalues, that is called the {\it eigenvalue geodesic polygon}. 
The structure theorem \ref{convex} and \ref{leaf-eigen} lead that the leaf of an endomorphism of dimension at least $3$ includes the eigenvalue geodesic polygon. 

\thmenv{mainthm}
{
[Normal case]
\label{normal-polygon}
Leaf of a normal endomorphism is the eigenvalue geodesic polygon. 
}

\thmenv{remark}
{
Although not covered in this paper, leaf can also be defined for any bounded linear operator on any infinite dimensional inner product space, by taking the closure, and the above structure theorems also hold. 
}

In \S 2-3, there are some basic facts on a normal endomorphism and on the Poincar\'e metric. 
In \S 4, we define the leaf of an endomorphism. 
In \S 5, we describe the difference between numerical range and leaf. 
We decide the leaf of an endomorphism of dimension $1$ or $2$ in \S 6, and prove the structure theorems for the leaf of any endomorphism of dimension at least $3$ in \S 7-8. 

\setcounter{mainthm}{0}

\section{Normal endomorphisms}

Let $V$ be a real inner product space of finite dimension, and $\varphi$ a normal endomorphism on $V$. 
Then $V$ is the orthogonal direct sum of $\varphi$-stable subspaces of dimension $1$ or $2$ (see Theorem 10.10 in \cite{R}). 
In this section, we provide a brief introduction to this fact. 

Put $V_\C=V\otimes_\R\C=V\oplus\unitCpx V$ the coefficient extension of $V$ to $\C$. 
There exists a Hermitian inner product on $V_\C$ which extends the real inner product on $V$. 
The same dot symbol is used for both inner products.
Each $\vc v\in V_\C=V\oplus\unitCpx V$ is expressed by $\vc v_r+\unitCpx\vc v_i$ for some $\vc v_r$, $\vc v_i\in V$. 
$\vc v_r$ is called the real part of $\vc v$, denoted by $\real\vc v$, and $\vc v_i$ is called the imaginary part of $\vc v$, denoted by $\imag\vc v$. 
$\overline{\vc v}=\vc v_r-\unitCpx\vc v_i$ is called the complex conjugate of $\vc v$. 
Each endomorphism $\varphi$ on $V$ can be extended to the endomorphism on $V_\C$ over $\C$. 
The same symbol $\varphi$ is used for both endomorphisms.

Let $\vc v\in V_\C$ be an eigenvector of $\varphi$ corresponding to an eigenvalue $\lambda\in\C$. 
\begin{align*}
&\varphi(\vc v)=\varphi(\real\vc v+\unitCpx\,\imag\vc v)=\varphi(\real\vc v)+\unitCpx\,\varphi(\imag\vc v) \\
&\varphi(\vc v)=\lambda\,\vc v
=(\real\lambda\,\real\vc v-\imag\lambda\,\imag\vc v)+\unitCpx\,(\imag\lambda\,\real\vc v+\real\lambda\,\imag\vc v) \\
&\quad \therefore\ \ (\varphi(\real\vc v),\ \varphi(\imag\vc v))=(\real\vc v,\ \imag\vc v)\,(\real\lambda\,\mtE+\imag\lambda\,\mtJ)
\end{align*}
where $\mtE$ and $\mtJ$ are the quadratic matrices defined in Example \ref{eg-leaf}. 
The complex conjugate $\overline{\vc v}$ is an eigenvector of $\varphi$ corresponding to the eigenvalue $\overline{\lambda}$. 
If $\lambda\in\R$, then $\real\vc v$, $\imag\vc v$ and $\overline{\vc v}$ are eigenvectors of $\varphi$ corresponding to $\lambda$. 
If $\lambda\not\in\R$, 
since $\vc v$ and $\overline{\vc v}$ are linearly independent on $\C$, 
$\{\real\vc v,$ $\imag\vc v\}$ is a basis of the $\varphi$-stable subspace $\langle\real\vc v,\,\imag\vc v\rangle=(\C\,\vc v+\C\,\overline{\vc v})\cap V$. 
Especially, if $\vc v$ and $\overline{\vc v}$ are orthogonal in $V_\C$, then $\real\vc v$ and $\imag\vc v$ are orthogonal in $V$. 
Because 
\[\real \vc v\cdot\imag \vc v=\frac{\vc v+\overline{\vc v}}{2}\cdot\frac{\vc v-\overline{\vc v}}{2\unitCpx}=\frac{-1}{4}(\|\vc v\|^2-\vc v\cdot\overline{\vc v}+\overline{\vc v}\cdot\vc v-\|\overline{\vc v}\|^2)=\frac{1}{2}\real(\vc v\cdot\overline{\vc v})=0. \]
$\{\real \vc v/\|\real\vc v\|,\,\imag \vc v/\|\imag\vc v\|\}$ is an orthonormal basis of $\langle\real\vc v,\,\imag\vc v\rangle$. 
Moreover, since $\|\real\vc v\|$ $=\|\imag\vc v\|=\|\vc v\|/\ssqrt{2}$, 
the matrix $\real\lambda\,\mtE+\imag\lambda\,\mtJ$ is the matrix representation of the restriction of $\varphi$ to the subspace $\langle\real\vc v,\,\imag\vc v\rangle$ with respect to the orthonormal basis $\{\real\vc v/\|\real\vc v\|,\,\imag\vc v/\|\imag\vc v\|\}$. 

The adjoint endomorphism $\varphi^\ast$ of $\varphi$ is defined by the endomorphism on $V$ satisfying that $\varphi(\vc v_1)\cdot\vc v_2=\vc v_1\cdot\varphi^\ast(\vc v_2)$ for any $\vc v_1,\vc v_2\in V$. 
The adjoint endomorphism of $\varphi$ as an endomorphism on $V_\C$ is equal to the extension to $V_\C$ of $\varphi^\ast$ on $V$. 
An endomorphism that commutative with its adjoint is said to be normal. 
Each normal endomorphism on $V_\C$ is diagonalizable by an orthonormal basis consisting eigenvectors on $V_\C$. 
Every normal endomorphism on $V$ has as the matrix representation the direct sum of the real eigenvalues and of the $a\,\mtE{+}b\,\mtJ$-type quadratic matrices ($a,b\in\R$, $b\neq0$) corresponding to the non-real eigenvalues ($a+b\unitCpx$), given an appropriate orthonormal basis. 

\section{The Poincar\'e metric on \boldmath{$\Hu$}}

The Poincar\'e metric on the complex upper half plane $\Hu$ is defined by $ds=|dz|/\imag z$ for $z\in\Hu$. 
The distance $\rho(z_1,\,z_2)$ between two points $z_1$, $z_2\in\overline{\Hu}$ is given by 
\[\rho(z_1,\,z_2)=\int_{z_1}^{z_2}|ds|=\log\frac{1+\delta(z_1,\,z_2)}{1-\delta(z_1,\,z_2)} \qquad\text{where}\quad \delta(z_1,\,z_2)=\frac{|z_1-z_2|}{|z_1-\overline{z_2}|}.\] 
We see that $0\le\delta(z_1,\,z_2)\le1$, and the function $\log\frac{1+k}{1-k}$ is monotonic increasing on $k$. 
$\R\cup\{\infty\}$ is called the line at infinity, because $\delta(z_1,\,z_2)=1$ and $\rho(z_1,\,z_2)=\infty$, if and only if at least one of $z_1$ and $z_2$ lies on $\R\cup\{\infty\}$. 
The geodesic on $\Hu$ with respect to the Poincar\'e metric is a part of a line parallel to the imaginary axis or a part of a semicircle whose center is on the real axis. 
A geodesic circle whose center is $z_0\in\Hu$ and radius is $r=\log\frac{1+k}{1-k}$ for $0\le k\le1$ is defined by 
\[\gC(z_0;\,k)=\{z\in\Hu\mid\delta(z,\,z_0)=k\}=\{z\in\Hu\mid\rho(z,\,z_0)=r\}.\]
By putting $z=x+y\unitCpx$, $z_0=p+q\unitCpx\in\Hu\cup\R$ ($x,y,p,q\in\R$, $y\ge0$, $q\ge0$), the equation $\delta(z,\,z_0)=k$ or $\rho(z,\,z_0)=r$ is equivalent to 
\[(x-p)^2+(y-\frac{1+k^2}{1-k^2}q)^2=((\frac{1+k^2}{1-k^2})^2-1)\,q^2=(\frac{2\,k}{1-k^2}\,q)^2. \]
Therefore, the geodesic circle $\gC(z_0;\,k)$ is a ordinary circle in $\C$ with center $p+\frac{1+k^2}{1-k^2}q\unitCpx$ and radius $\frac{2\,k}{1-k^2}q$. 
A geodesic circle on the Poincar\'e metric is simply called a geodesic circle, and an ordinary circle in the complex plane is simply called a circle. 

\section{Leaf and real eigenvalues}

Let $\varphi$ be an endomorphism on a real inner product space $V$ of finite dimension. 
The continuous map $\pv$ on $V^\times$ ($=V\smallsetminus\{\vc 0\}$) to $\C$ defining by 
\[\pv(\vc v)=\frac{\|\varphi(\vc v)\|}{\|\vc v\|}\,e^{\sqrt{-1}\,\theta(\vc v,\,\varphi(\vc v))}\qquad (\vc v\in V^\times),\] 
where the angle $\theta=\theta(\vc v,\,\varphi(\vc v))\in[0,\,\pi]\subset\R$ between $\vc v$ and $\varphi(\vc v)$ is gdefined by $\vc v\cdot\varphi(\vc v)=\|\vc v\|\,\|\varphi(\vc v)\|\,\cos\theta$. 
When $\varphi(\vc v)=\vc 0$, we don't define the angle $\theta$,but can put $\pv(\vc v)=0$ because of $\|\varphi(\vc v)\|=0$. 
$\pv(k\vc v)=\pv(\vc v)$ holds for any $k\in\R^\times$, because $\|\varphi(k\vc v)\|/\|k\vc v\|=\|\varphi(\vc v)\|/\|\vc v\|$ and $\theta(k\vc v,\,\varphi(k\vc v))=\theta(\vc v,\,\varphi(\vc v))$. 
So the continuous map $\pv$ is considered as the continuous map on the projective space $\Proj(V)=V^\times/\R^\times$. 

For any $\vc v_1$, $\vc v_2\in V^\times$, put $\sigma(\vc v_1,\,\vc v_2)=\sqrt{\|\vc v_1\|^2\,\|\vc v_2\|^2-(\vc v_1\cdot\vc v_2)^2}$. 
Then we have 
\[\pv(\vc v)=\frac{\vc v\cdot\varphi(\vc v)+\sigma(\vc v,\,\varphi(\vc v))\,\unitCpx}{\|\vc v\|^2}.\] 
The real part of $\pv(\vc v)$ is equal to the Rayleigh quotient $\frac{\vc v\cdot\varphi(\vc v)}{\|\vc v\|^2}$ of $\varphi$, and the square of absolute value of $\pv(\vc v)$ to that of $\varphi^\ast\varphi$. 
The operator norm $\|\varphi\|$ of $\varphi$ is the supremum of the absolute value of $\pv(\vc v)$. 

\thmenv{definition}
{[leaf]
\label{def-leaf}
The image of $\pv$ 
\[\Psi(\varphi)=\pv(V^\times)=\pv(\Proj(V))\]
is called the {\it leaf} of the endomorphism $\varphi$. 
}

\thmenv{proposition}
{\label{leaf-arc}
$\Psi(\varphi)$ is path connected. 
}

\thmenv{proof}
{
Let $\lambda_1$, $\lambda_2\in\Psi(\varphi)$, and take $\vc v_1$, $\vc v_2\in V^\times$ satisfying that $\lambda_1=\pv(\vc v_1)$ and $\lambda_2=\pv(\vc v_2)$. 
The path $[0,\,1]\ni t\mapsto \pv((1-t)\,\vc v_1+t\,\vc v_2)\in\Psi(\varphi)$ in $\Psi(\varphi)$ connects $\lambda_1$ and $\lambda_2$. 
Hence, the leaf $\Psi(\varphi)$ is path connected. 
}

The leaf $\Psi(\varphi)$ is a bounded closed set in $\C$ because the projective space $\Proj(V)$ is compact. 
Since an argument $\theta$ belongs to the closed interval $[0,\,\pi]$, the value of $\pv$ belongs to the closure of the upper half plane $\overline{\Hu}=\Hu\cup\R\cup\{\infty\}$. 
The leaf is included in $\overline{\Hu}$. 
We will see that the Poincar\'e metric of $\Hu$ plays an important role on the leaf in the later chapters.
In closing this section, we show the following theorem on the real eigenvalues. 

\thmenv{theorem}
{\label{leaf-real-eigen}
The intersection of the leaf and the real axis is the set of real eigenvalues. 
}

\thmenv{proof}
{
Let $\varphi$ an endomorphism on a real inner product space of finite dimension $V$. 
Let $\lambda\in\R$ a real eigenvalue of $\varphi$, and $\vc v\in V$ an eigenvector of $\varphi$ corresponding to $\lambda$. 
If $\lambda=0$, then $\varphi(\vc v)=\vc 0$, so $\pv(\vc v)=0=\lambda$. 
Assume that $\lambda\neq0$. 
Since $\varphi(\vc v)=\lambda\,\vc v$ is parallel to $\vc v$, the angle $\theta=\theta(\vc v, \varphi(\vc v))$ is equal to $0$ if $\lambda>0$, and to $\pi$ if $\lambda<0$. 
Since $\|\varphi(\vc v)\|/\|\vc v\|=\|\lambda\,\vc v\|/\|\vc v\|=|\lambda|$, and $e^{\theta\sqrt{-1}}={\rm sign}\lambda=\lambda/|\lambda|$, we have that $\pv(\vc v)=\lambda$. 
Thus $\lambda=\pv(\vc v)\in\Psi(\varphi)\cap\R$. 

Conversely, we will show that each element $\lambda\in\Psi(\varphi)\cap\R$ is an eigenvalue of $\varphi$. 
Take $\vc v\in V^\times$ satisfying that $\lambda=\pv(\vc v)$. 
In the case that $\lambda=0$, it follows from $\pv(\vc v)=\lambda=0$ that $\|\varphi(\vc v)\|=0$, i.e. $\varphi(\vc v)=\vc 0=0\,\vc v$, 
so $0$ is an eigenvalue of $\varphi$. 
Assume that $\lambda\neq0$. 
Since $\|\varphi(\vc v)\|=|\pv(\vc v)|\,\|\vc v\|=|\lambda|\,\|\vc v\|\neq0$, the angle $\theta=\theta(\vc v,\,\varphi(\vc v))$ is defined. 
And $\theta$ is also the angle of the real number $\lambda=\pv(\varphi)$, so $\theta=0$ or $\pi$. 
$\varphi(\vc v)$ is parallel to $\vc v$, namely $\vc v$ is an eigenvector of $\varphi$. 
Put $\lambda'$ be the eigenvalue of $\varphi$ satisfying that $\varphi(\vc v)=\lambda'\vc v$, then it follows from the first step in this proof that $\pv(\vc v)=\lambda'$. 
Hence $\lambda=\pv(\vc v)=\lambda'$ is an eigenvalue of $\varphi$.  
}

\section{Numerical range vs. leaf}

Well-known object as a convex domain containing all eigenvalues is a numerical range (cf. \cite{ERS}, \cite{GR}). 
Here we introduce the fundamental theorems that were the starting point for studies of numerical range, according to chapter 1 of Guatafson-Rao’s text boot (\cite{GR}). 

Let $\varphi$ be a bounded linear operator on a Hermite inner product space $V$. 
The numerical range of $\varphi$ is the subset in $\C$ defined by 
\[W(\varphi)=\{\frac{\vc v\cdot\varphi(\vc v)}{\|\vc v\|^2}\mid\vc v\in V^\times\}=\{\vc v\cdot\varphi(\vc v)\mid \|\vc v\|=1\},\] 
that is the range of the Rayleigh quotients $\vc v\cdot\varphi(\vc v)/\|\vc v\|^2$, $\vc v\in V^\times$. 
The fundamental theorems of numerical range are following. 

\thmenv{theorem}
{[Convexity, Toeplitz-Hausdorff, \cite{H}, \cite{T}]
\label{Toe-Hau}
The numerical range of an operator is convex with respect to the Euclidean metric on $\C$. 
}
\thmenv{theorem}
{[Spectral inclusion, \cite{W}]
\label{Win}
The spectrum of a bounded linear operator is contained in the closure of its numerical range.
}
\thmenv{theorem}
{[\cite{B}, \cite{S}] 
\label{polygon-NR}
The closure of numerical range of a normal bounded linear operator is the convex full of its spectrum. 
}
\thmenv{theorem}
{[Equivalent norm]
\label{opNorm-NR}
Put $w(\varphi)=\sup\{|\lambda|\mid\lambda\in W(\varphi)\}$ (the numerical radius), and $\|\varphi\|=\sup\{\|\varphi(\vc v)\|/\|\vc v\|\mid\vc v\in V^\times\}$ (the operator norm). 
Then $w(\varphi)\le\|\varphi\|\le2w(\varphi)$. 
}

The structure theorems of leaf (Theorem \ref{convex}, \ref{leaf-eigen}, \ref{normal-polygon}) correspond to the first three of the fundamental theorems above.
The corresponding result to the last theorem on the leaf is that 
\[\|\varphi\|=\sup\{|\lambda|\mid\lambda\in\varPsi(\varphi)\}\]
which follows from the definition of leaf. 

There is an essential difference between Theorem \ref{Toe-Hau} and Theorem \ref{convex}: numerical range is convex with respect to the Euclidean metric, and leaf is convex with respect to the Poincar\'e metric. 
Different metrics lead to different appearances of the real eigenvalues in the figures. 
The intersection of numerical range and the real axis is an interval containing all real eigenvalues (Theorem \ref{Toe-Hau}, \ref{Win}), whereas that of leaf is the set of all real eigenvalues (Theorem 4.3). 
The information on each of the real eigenvalues is lost in numerical range, but retained in leaf.

Theorem \ref{leaf-eigen} and Theorem \ref{Win} have different ways of proving. 
In the case of numerical range, the Rayleigh quotient of an eigenvector is the eigenvalue. 
But, since any endomorphism on a real vector space has no real eigenvector with respect to the non-real eigenvalue which is a non-real root of the characteristic polynomial, the same method cannot be use for leaf as for numerical range. 
We need to use Theorem \ref{convex}, to find a real vector for which the value of the function $\pv$ is the non-real eigenvalue (Proposition \ref{leaf-simp-eigen}). 

Theorem \ref{convex} and Theorem \ref{Toe-Hau} also have different ways to proving. 
The line connecting the two points of numerical range is given by the image by the Rayleigh quotient of real coefficient linear combination of two vectors giving the two points on numerical range. 
The geodesic connecting the two points of leaf is contained in the interior of the closed curve formed by the image by $\pv$ of the subspace generated by the two vectors giving the two points (Proposition \ref{loop}). 
It can be shown that the geodesic connecting the two points is contained in leaf, by continuously moving that closed curve using a third vector to fill the interior of that closed curve (Proposition \ref{loop-filled}). 

In the case of a normal endomorphism, both numerical range and leaf are convex sets containing all eigenvalues. 
The metrics are different, so the geodesics are different shapes. 
The edges which are geodesics connecting the eigenvalues have different shapes for numerical range and leaf.

\section{\boldmath Leaf of an endomorphism of dimension $1$ or $2$}

In this section, we will show the following theorem.

\thmenv{theorem}
{\label{leaf-dim12}
Leaf of an endomorphism of dimension $1$ or $2$ is either a single point set, a circle or a bended circle. 
}

A bended circle is defined by a figure formed by bending a circle that intersects the real axis upward along the real axis. 
The center and the radius of a bended circle are defined by the center and the radius of the original circle, respectively. 
Let $C$ be a circle in $\C$ crossing the real axis, and put $C'$ be the complex conjugate of $C$. 
Then $(C\cup C')\cap\overline{\Hu}$ is a bended circle. 

Clearly leaf of each endomorphism of dimension $1$ is a single point set. 

\thmenv{proposition}
{\label{leaf-dim2}
Let $V$ be a real inner product space of dimension 2, and $\varphi$ an endomorphism on $V$. 
\begin{enumerate}[ $($\rm a$)$]
\item 
If $\varphi$ has $\matrixform{cc}{a&b\\c&d}$ as a matrix representation with respect to an orthonormal basis of $V$, then the leaf $\Psi(\varphi)$ contains $a+|b|\unitCpx$, $a+|c|\unitCpx$, $d+|b|\unitCpx$ and $d+|c|\unitCpx$. 
\item 
Put $\alpha=\frac{a+d}{2}+\frac{|b-c|}{2}\unitCpx$ and $r=\frac{1}{2}\sqrt{(a-d)^2+(b+c)^2}$. 
Both $\alpha$ and $r$ are independent of choice of orthonormal basis of $V$. 
\item 
Let $C\subset\C$ be the circle with center $\alpha$ and radius $r$. 
If $\varphi$ does not have any real eigenvalues, then $\Psi(\varphi)=C$. 
Moreover, if $r>0$, then the eigenvalue with positive imaginary part is located on the circumference or on the interior of $C$. 
\item 
If all eigenvalues of $\varphi$ are real, then $C$ crosses the real axis at the eigenvalues, and $\Psi(\varphi)$ is a bended circle $(C\cup C')\cap\overline{\Hu}$, 
where $C'$ is the complex conjugate of $C$. 
\end{enumerate}
}

\thmenv{proof}
{
\thmNum Let $\{\vc u_1,\,\vc u_2\}$ be an orthonormal basis of $V$ satisfying that $\matrixform{cc}{a&b\\c&d}$ is the matrix representation of $\varphi$ with respect to this basis. 
Then $\varphi(\vc u_1)=a\,\vc u_1+c\,\vc u_2$. 
Since $\|\vc u_1\|=1$, $\|\varphi(\vc u_1)\|^2=a^2+c^2$, $\vc u_1\cdot\varphi(\vc u_1)=a$ and $\sigma(\vc u_1,\,\varphi(\vc u_1))=|c|$, we have $\pv(\vc u_1)=a+|c|\unitCpx$. 
Similarly, we have $\pv(\vc u_2)=d+|b|\unitCpx$, $\pv((a-d)\vc u_1+(b+c)\vc u_2)=a+|b|\unitCpx$ and $\pv((b+c)\vc u_1-(a-d)\vc u_2)=d+|c|\unitCpx$. 
Their four points belongs to $\Psi(\varphi)$. 
\\
\thmNum Let $A=\matrixform{cc}{a&b\\c&d}$ and $A'=\matrixform{cc}{a'&b'\\c'&d'}$ be the matrix representations of $\varphi$ with respect to the two orthonormal bases of $V$. 
Put $P$ be the translation matrix between their orthonormal bases, then $A'=P^{-1}AP=\tp PAP$, because $P$ is an orthogonal matrix. 
The column vectors $\vc p$ and $\vc q$ where $P=(\vc p\ \ \vc q)$ are an orthonormal basis of $\R^2$. 
Comparing the $(1,2)$ and $(2,1)$ components of $A'=\tp PAP$ yields $b'=\tp \vc pA\vc q$ and $c'=\tp \vc qA\vc p=\tp \vc p\tpp A\vc q$. 
And then $b'-c'=\tp \vc p(A-\tpp A)\vc q=\tp \vc p(b-c)\mtJ\vc q=(b-c)\tp \vc p\mtJ\vc q$, where $\mtJ$ is defined in Example \ref{eg-leaf}. 
$\mtJ$ is the $\pi/2$-rotation on $\R^2$, so $\mtJ\vc q=\pm\vc p$. 
Hence, $|b'-c'|=|(b-c)(\pm\|\vc p\|^2)|=|b-c|$. 

The trace and the determinant of a matrix representation are independent of choice of basis.
That is $a+d=a'+d'$ and $ad-bc=a'd'-b'c'$, then 
\[\alpha=\frac{a+d}{2}+\frac{|b-c|}{2}\unitCpx=\frac{a'+d'}{2}+\frac{|b'-c'|}{2}\unitCpx \]
and
\begin{align*}
(2r)^2 & =(a-d)^2+(b+c)^2=(a+d)^2+(b-c)^2-4(ad-bc)\\
& =(a'+d')^2+(b'-c')^2-4(a'd'-b'c')=(a'-d')^2+(b'+c')^2.
\end{align*} 
We have that $\alpha$ and $r$ are independent of choice of orthonormal basis of $V$. 
\\
\thmNum Let $D$ be the discriminant of the characteristic polynomial of $\varphi$. 
$\varphi$ does not have any real eigenvalues, so $D$ is negative.  
Since $D=(a+d)^2-4(ad-bc)=(a-d)^2+4bc$, we have $bc<0$. 
Four points of (a) form the rectangle with center $\alpha$ and diagonal length $2r$. 
Therefore, four points in (a) lies on $C$. 

Let $\lambda'=\pv(\vc v)\in\Psi(\varphi)$ where $\vc v\in V^\times$, and put $\vc u=\vc v/\|\vc v\|$.  
There exists an orthonormal basis of $V$ containing $\vc u$. 
Since $C$ is independent of choice of orthonormal basis of $V$, $\lambda'=\pv(\vc v)=\pv(\vc u)$ lies on $C$. 
$\Psi(\varphi)$ is included in $C$. 

We will show that $C\subset\Psi(\varphi)$. 
The leaf is path connected (Proposition \ref{leaf-arc}), so the complement $C\smallsetminus\Psi(\varphi)$ is at most $1$ connected component. 
For a matrix representation $\matrixform{cc}{a'&b'\\c'&d'}$ of $\varphi$ with respect to an orthonormal basis of $V$, two points $a'+|c'|\unitCpx$ and $d'+|b'|\unitCpx\in\Psi(\varphi)$ are on the diagonal on $C$. 
$\Psi(\varphi)$ is point symmetric with respect to the center $\alpha$ of $C$. 
The complement $C\smallsetminus\Psi(\varphi)$ is also point symmetric. 
Thus, the complement $C\smallsetminus\Psi(\varphi)$ is the empty set or whole of $C$. 
Since $\Psi(\varphi)$ is not empty, $\C\subset\Psi(\varphi)$. 
So $\Psi(\varphi)=C$. 

Let $\lambda$ be an eigenvalue of $\varphi$ with positive imaginary part. 
Then we have $\lambda=\frac{a+d}{2}+\frac{\sqrt{D}}{2}$. 
Since $D+(b-c)^2=(a-d)^2+(b+c)^2=(2r)^2\ge0$ and $D<0$, we have that $\sqrt{|D|}\le|b-c|$. 
The distance between $\lambda$ and $\alpha$ (the center of $C$) is at most the radius of $C$, as below. 
\[|\lambda-\alpha|^2=\frac{(\sqrt{|D|}-|b-c|)^2}{4}=\frac{|D|+(b-c)^2-2\sqrt{|D|}\,|b-c|}{4}\le\frac{|D|+(b-c)^2-2|D|}{4}=r^2.\]
Namely, $\lambda$ is located on the circumference or on the interior of $C$. 
\\
\thmNum 
Let $\lambda$ be a real eigenvalue of $\varphi$, which is denoted by $\frac{a+d\pm\sqrt{D}}{2}$. 
\[|\lambda-\alpha|^2 =(\lambda-\frac{a+d}{2})^2+(\frac{|b-c|}{2})^2=\frac{D+(b-c)^2}{4}=r^2. \]
Hence, each real eigenvalue lies on the circumference of the circle $C$. 
Let $C'$ be the circle with center $\overline{\alpha}=\frac{a+d}{2}-\frac{b-c}{2}\unitCpx$ and radius $r$. 
$C'$ is symmetric with $C$ about the real axis.
Hence, each real eigenvalue also lies on the circumference of the circle $C'$. 
We consider eight points $a\pm b\unitCpx$, $a\pm c\unitCpx$, $d\pm b\unitCpx$ and $d\pm c\unitCpx$. 
Four points $a+b\unitCpx$, $a-c\unitCpx$, $d+b\unitCpx$ and $d-c\unitCpx$ are the vertices of the rectangle with center $\alpha$ and diagonal length $2r$, and are on $C$. 
Another four points $a-b\unitCpx$, $a+c\unitCpx$, $d-b\unitCpx$ and $d+c\unitCpx$ are the vertices of the rectangle with center $\overline{\alpha}$ and diagonal length $2r$, and are on $C'$. 
Therefore these eight points are on the union of the circles $C$ and $C'$. 
Four of eight points, $a+|b|\unitCpx$, $a+|c|\unitCpx$, $d+|b|\unitCpx$ and $d+|c|\unitCpx$ belongs to $\overline{\Hu}$, and then are on the bended circle $(C\cup C')\cap\overline{\Hu}$ which is the shape made by bending the circle $C$ upward along the real axis. 
We have that $\Psi(\varphi)$ is included in $(C\cup C')\cap\overline{\Hu}$. 
The reverse inclusive relation $(C\cup C')\cap\overline{\Hu}\subset\Psi(\varphi)$ is valid in the same way as the proof of (c). 
}

\thmenv{proposition}
{\label{subsp-dim2-real}
Let $\varphi$ be an endomorphism on a real inner product space $V$ of dimension at least $2$. 
Assume that $\varphi$ has two different real eigenvalues $\lambda_1$ and $\lambda_2$. 
Let $\vc v_1$ and $\vc v_2\in V$ be unit eigenvectors of $\varphi$ corresponding to $\lambda_1$ and $\lambda_2$, respectively, and put $\beta=\vc v_1\cdot\vc v_2$. 
Then $\pv(\langle\vc v_1,\,\vc v_2\rangle^\times)$ is the bended circle with center $\alpha$ and radius $r$, where $\alpha=\frac{\lambda_1+\lambda_2}{2}+t\unitCpx$, and $t$ and $r$ are nonnegative real numbers satisfying that 
\[r^2=(\frac{\lambda_1-\lambda_2}{2})^2+t^2, \qquad t=\frac{|\lambda_1-\lambda_2|}{2}\,\frac{|\beta|}{\sqrt{1-\beta^2}}, \quad\text{and}\quad r=\frac{|\lambda_1-\lambda_2|}{2}\,\frac{1}{\sqrt{1-\beta^2}}.\]
}

\thmenv{proof}
{
Since the last proposition leads that $\lambda_1$ and $\lambda_2$ lies on the circle with center $\alpha$ and radius $r$, we have that $|\lambda_1-\alpha|=|\lambda_2-\alpha|=r$, namely 
\[(\frac{\lambda_1-\lambda_2}{2})^2+t^2=r^2\]
Put $\vc w=\vc v_1+\vc v_2\in \langle\vc v_1,\,\vc v_2\rangle^\times$, then 
\begin{gather*}
\|\vc w\|^2=2(1+\beta),\quad 
\|\varphi(\vc w)\|^2=\lambda_1^2+\lambda_2^2+2\lambda_1\lambda_2\beta, \quad \text{and} \quad 
\vc w\cdot\varphi(\vc w)=(\lambda_1+\lambda_2)(1+\beta) \\
\therefore \, |\pv(\vc w)|^2=\frac{\|\varphi(\vc w)\|^2}{\|\vc w\|^2}=\frac{\lambda_1^2+\lambda_2^2+2\lambda_1\lambda_2\beta}{2(1+\beta)}, \quad 
\real\pv(\vc w)=\frac{\vc w\cdot\varphi(\vc w)}{\|\vc w\|^2}=\frac{\lambda_1+\lambda_2}{2}
\end{gather*}
Hence 
\[
(\imag\pv(\vc w))^2
=|\pv(\vc w)|^2-(\real\pv(\vc w))^2
=(\frac{\lambda_1-\lambda_2}{2})^2\,\frac{1-\beta}{1+\beta}. 
\]
Since the real part of $\pv(\vc w)$ coincides with that of the center $\alpha$ of the bended circle $\pv(\langle\vc v_1,\,\vc v_2\rangle^\times)$. 
\[\imag\pv(\vc w)=\text{(the radius of bended circle)}\pm\imag\alpha=r\pm t.\]
We obtain that 
\[r^2-t^2=(\frac{\lambda_1-\lambda_2}{2})^2, \quad \text{and} \quad r\pm t=\frac{\lambda_1-\lambda_2}{2}\sqrt{\frac{1-\beta}{1+\beta}}\]
Therefore
\begin{gather*} 
r=\frac{\lambda_1-\lambda_2}{2}\,\frac{1}{2}\bigg(\sqrt{\frac{1-\beta}{1+\beta}}+\sqrt{\frac{1+\beta}{1-\beta}}\bigg)=\frac{\lambda_1-\lambda_2}{2}\frac{1}{\sqrt{1-\beta^2}},\\
 t=\frac{\lambda_1-\lambda_2}{2}\,\frac{1}{2}\,\bigg|\,\sqrt{\frac{1-\beta}{1+\beta}}-\sqrt{\frac{1+\beta}{1-\beta}}\ \bigg|\,=\frac{\lambda_1-\lambda_2}{2}\frac{|\beta|}{\sqrt{1-\beta^2}}.
\end{gather*}
}

\section{Leaf in \boldmath{$\overline{\Hu}$} and eigenvalues with nonnegative imaginary part}

There are various shapes within the leaves of endomorphisms of dimension at least $3$. 
We will show the structure theorems on leaf. 
The first one 
is the geometric property that each leaf is convex on the Poincar\'e metric, and then each leaf is simply connected. 
The second one 
is the algebraic property that each leaf contains all eigenvalues with nonnegative complex part. 

Let $V$ be a real inner product space of finite dimension, and $\varphi$ an endomorphism on $V$. 
Let $W$ be a $\varphi$-stable subspace in $V$. 
Since the restriction $\varphi|_W$ of $\varphi$ on $W$ is the endomorphism on $W$, the subset $\pv(W^\times)$ of $\Psi(\varphi)$ is equal to the leaf $\Psi(\varphi|_W)$ of $\varphi|_W$. 
$\pv(W^\times)=\Psi(\varphi|_W)$ is called the subleaf corresponding to $W$ and denoted by $\Psi(\varphi;\,W)$. 

\thmenv{proposition}
{\label{leaf-subsp-dim2}
A subleaf corresponding to a $2$-dimensional stable subspace is one of three kinds of figures on Theorem \ref{leaf-dim12}.
}

\thmenv{proposition}
{\label{leaf-simp-eigen}
If the leaf is simply connected, then it contains all eigenvalues with nonnegative imaginary part. 
}

\thmenv{proof}
{
It is already shown that the leaf contains all real eigenvalues, by Theorem \ref{leaf-real-eigen}. 
Let $\lambda\in\Hu$ be an eigenvalue of $\varphi$ with positive imaginary part, and $\vc v\in V_\C$ an eigenvector of $\varphi$ corresponding to $\lambda$. 
$W=\langle\real\vc v,\,\imag\vc v\rangle$ is a $\varphi$-stable subspace of dimension $2$. 
By Propositions \ref{leaf-subsp-dim2} and \ref{leaf-dim2}, $\lambda$ is located on the circumference or on the inside of the circle $\Psi(\lambda;\,W)=\pv(W^\times)$. 
By the assumption that $\Psi(\varphi)$ is simply connected, $\Psi(\varphi)$ includes the circumference or on the interior of the circle $\Psi(\lambda;\,W)$. 
Hence, $\lambda$ belongs to $\Psi(\varphi)$. 
}

Put $\ell(W)=\pv(W^\times)=\pv(\Proj(W))$ for any $2$-dimensional subspace $W$ in $V$. 
$\ell(W)$ is a closed curve, because the projective space $\Proj(W)=W^\times/\R^\times$ is homeomorphic to the circle. 
Let $\D(W)$ be the simple connected domain bounded by the closed curve $\ell(W)$, 
and put $\D(W)^\circ=\D(W){\smallsetminus}\ell(W)$ the interior of $\D(W)$. 

\thmenv{proposition}
{\label{loop}
Let $\lambda_1$, $\lambda_2\in\Psi(\varphi)\subset\overline{\Hu}$ with $\lambda_1\neq\lambda_2$. 
Take $\vc v_1$, $\vc v_2\in V^\times$ with $\|\vc v_1\|=\|\vc v_2\|=1$, $\lambda_1=\pv(\vc v_1)$ and $\lambda_2=\pv(\vc v_2)$, and put $W=\langle\vc v_1,\vc v_2\rangle$ which is a $2$-dimensional subspace of $V$. 
Let $s(\lambda_1,\,\lambda_2)$ be the geodesic connecting $\lambda_1$ to $\lambda_2$ on $\overline{\Hu}$. 
\begin{enumerate}[ $($\rm a$)$]
\item 
The geodesic $s(\lambda_1,\lambda_2)$ is included in the domain $\D(W)$. 
\item 
If $s(\lambda_1,\,\lambda_2)$ does not pass through the interior $\D(W)^\circ(=\D(W){\smallsetminus}\ell(W))$, then the interior $\D(W)^\circ$ is empty, and $\ell(W)(=\partial\D(W))$ is a geodesic including $s(\lambda_1,\,\lambda_2)$. 
\item 
If $s(\lambda_1,\,\lambda_2)$ passes through the interior $\D(W)^\circ$, then $s(\lambda_1,\,\lambda_2)\cap\ell(W)=\{\lambda_1,\,\lambda_2\}$. 
\item 
If the interior $\D(W)^\circ$ is not empty, then $\ell(W)$ is a closed simple curve. 
\end{enumerate}
}

\thmenv{proof}
{
We observe the positional relationship between the closed curve $\ell(W)$ and the geodesic $s(\lambda_1,\,\lambda_2)$. 
Express a typical $\vc w\in\Proj(W)=W^\times/\R^\times$ is the form by $x\vc v_1+\vc v_2$ for some $x\in\Proj^2(\R)=(\R^2)^\times/\R^\times=\R\cup\{\infty\}$. 
Note that $\vc v_1=\infty\vc v_1+\vc v_2$ in this notation. 
We have that 
\begin{gather*}
\|\vc w\|^2=x^2+2x(\vc v_1\cdot\vc v_2)+1, \\
\|\varphi(\vc w)\|^2=x^2\|\varphi(\vc v_1)\|^2+2x(\varphi(\vc v_1)\cdot\varphi(\vc v_2))+\|\varphi(\vc v_2)\|^2, \\
\vc w\cdot\varphi(\vc w)=x^2(\vc v_1\cdot\varphi(\vc v_1))+x(\vc v_1\cdot\varphi(\vc v_2)+\vc v_2\cdot\varphi(\vc v_1))+\vc v_2\cdot\varphi(\vc v_2). 
\end{gather*}

A geodesic on the Poincar\'e metric is a part of a line parallel to the imaginary axis and semicircle with center on the real axis. 
$s(\lambda_1,\,\lambda_2)$ is a part of a line, if $\real\lambda_1=\real\lambda_2$, or a part of a semicircle, otherwise. 

(i) In the case that $\real\lambda_1=\real\lambda_2$. 
Put $a=\real\lambda_1=\real\lambda_2$ ($=\vc v_1\cdot\varphi(\vc v_1)=\vc v_2\cdot\varphi(\vc v_2)$). 
$s(\lambda_1,\,\lambda_2)$ is a part of the line whose real part is $a$. 
\[\real\pv(\vc w)-a = \frac{\vc w\cdot\varphi(\vc w)}{\|\vc w\|^2}-\vc v_1\cdot\varphi(\vc v_1)
=R\,\frac{x}{\|\vc w\|^2}, \]
where 
\[R=\vc v_1\cdot\varphi(\vc v_2)+\vc v_2\cdot\varphi(\vc v_1)-a\,(\vc v_1\cdot\vc v_2). \]
Since $\|\vc w\|^2=O(x^2)$ as $x\to\infty$, we ramark that $\real\pv(x\,\vc v_1+\vc v_2)-a=O(x^{-1})$ as $x\to\infty$. 

(ii) In the case that $\real\lambda_1\neq\real\lambda_2$. 
$s(\lambda_1,\lambda_2)$ is a part of semicircle. 
Let $c\in\R$ be the center of the semicircle and $d$ the radius. 
Since 
\[d^2=|\lambda_i-c|^2=(\real\lambda_i-c)^2+(\imag\lambda_i)^2=|\lambda_i|^2-2c\,\real\lambda_i+c^2 \quad \text{($i=1,\,2$),}\]
we have that 
\[d^2-c^2=|\lambda_i|^2-2c\,\real\lambda_i, \qquad c=\frac{|\lambda_1|^2-|\lambda_2|^2}{2(\real \lambda_1-\real \lambda_2)}.\]
Thus, 
\[|\pv(\vc w)-c|^2-d^2=|\pv(\vc w)|^2-2c\,\real\pv(\vc w)
=R\,\frac{x}{\|\vc w\|^2}, \]
where 
\[R=2(\varphi(\vc v_1)\cdot\varphi(\vc v_2)+c(\vc v_1\cdot\varphi(\vc v_2)+\vc v_2\cdot\varphi(\vc v_1))+(c^2-d^2)(\vc v_1\cdot\vc v_2)). \]
In this case, we also remark that $|\pv(x\,\vc v_1+\vc v_2)-c|^2-d^2=O(x^{-1})$ as $x\to\infty$. 

In both cases (i) and (ii), the sign of $R\,t/\|\vc w\|^2$ shows the position on $\pv(\vc w)$ relative to $s(\lambda_1,\,\lambda_2)$. 
Let $K$ be the geodesic extending from the both endpoints of $s(\lambda_1,\,\lambda_2)$ to $\R\cup\{\infty\}$. 
$K$ divides the upper half plane $\overline{\Hu}\smallsetminus K$ into two connected areas $K_1$ and $K_2$. 
In the case that $K$ is a line, put $K_1$ be the right area of $K$, and $K_2$ the left area. 
In the case that $K$ is a semicircle, put $K_1$ be the outer area of $K$, and $K_2$ the inner area. 
Then $\pv(\vc w)$ belongs to $K$ if $R\,x/\|\vc w\|^2=0$, to $K_1$ if $R\,x/\|\vc w\|^2>0$, and to $K_2$ if $R\,x/\|\vc w\|^2<0$. 
Since $R$ is independent of $t$, and $\|\vc w\|$ is always positive, the sign change of $R\,t/\|\vc w\|^2$ coincides with the sign change of $t$. 
We separate two cases $R=0$ and $R\neq0$. 

(A) In the case that $R=0$, 
$\pv(\vc w)$ is always on $K$, so $\ell(W)$ is a part of $K$. 
$\ell(W)$ contains $\lambda_1$ and $\lambda_2$. 
Hence $s(\lambda_1,\,\lambda_2)\subset\ell(W)$. 

(B) In the case that $R\neq0$. 
If $x\neq0$, $\infty$, then $R\,x/\|\vc w\|^2\neq0$, so $\pv(\vc w)$ does not belong to $K$. 
$\pv(\vc w)$ belongs to one of $K_1$ or $K_2$ if $x>0$, and to the other if $x<0$. 
$\ell(W)$ goes around $s(\lambda_1,\,\lambda_2)$ as below: starts at $\lambda_1$ which is one of endpoints of the geodesic $s(\lambda_1,\,\lambda_2)$, goes through $K_1$, crosses $K$ at $\lambda_2$ which is the other endpoint of $s(\lambda_1,\,\lambda_2)$, goes through $K_2$, and returns to $\lambda_1$. 
Hence, $s(\lambda_1,\,\lambda_2)$ is included in $\D(W)$. 

In both cases (A) and (B), $s(\lambda_1,\,\lambda_2)$ is included in $\D(W)$, so (a) holds. 
(b) follows from (A), and (c) from (B). 

Finally, we will show that (d) holds. 
Let $\lambda_3$ be a self-intersection point of $\ell(W)$. 
If $\ell(W)$ has a normal crossing at $\lambda_3$, then there exist two geodesics that tangent to $\ell(W)$ at $\lambda_3$. 
Each geodesic between these geodesic intersects $\ell(W)$ in at least three points. 
This contradict (c), so there exists only one tangent geodesic at $\lambda_3$. 
There exists a geodesic that displaces the tangent geodesic just a little and crosses $\ell(W)$ by more than two points. 
This also contradict (c). 
Therefore, $\ell(W)$ does not intersect itself. 
}

\thmenv{proposition}
{\label{geod}
Let $\vc v_1$, $\vc v_2\in V^\times$. 
Assume that $\vc v_1$ and $\vc v_2$ satisfy the orthogonal condition: 
\[\langle\vc v_1,\,\varphi(\vc v_1)\rangle\perp\langle\vc v_2,\,\varphi(\vc v_2)\rangle.\]
Then the curve $\ell(\langle\vc v_1,\,\vc v_2\rangle)$ is the geodesic connecting $\pv(\vc v_1)$ to $\pv(\vc v_2)$. 
}

\thmenv{proof}
{
Put $\lambda_1=\pv(\vc v_1)$, $\lambda_2=\pv(\vc v_2)$ and $W=\langle\vc v_1,\,\vc v_2\rangle$. 
The orthogonal condition leads that $R=0$, where $R$ is defined in the proof of the last proposition. 
Hence, $\ell(W)$ is a geodesic connecting $\lambda_1$ to $\lambda_2$. 

Express a typical $\vc w\in V^\times$ is the form by $x\,\vc v_1+y\,\vc v_2$ for some $x$, $y\in\R$ with $(x,y)\neq(0,0)$. 
It follows from the orthogonal condition that 
\begin{gather*}
\|\vc w\|^2=x^2\|\vc v_1\|^2+y^2\|\vc v_2\|^2, \\
\|\varphi(\vc w)\|^2=x^2\|\varphi(\vc v_1)\|^2+y^2\|\varphi(\vc v_2)\|^2, \\
\vc w\cdot\varphi(\vc w)=x^2\vc v_1\cdot\varphi(\vc v_1)+y^2\vc v_2\cdot\varphi(\vc v_2). 
\end{gather*}
Hence, 
\begin{gather*}
\real\pv(\vc w)
=\frac{\vc w\cdot\varphi(\vc w)}{\|\vc w\|^2}
=\frac{(x\,\|\vc v_1\|)^2\,\real\lambda_1+(y\,\|\vc v_2\|)^2\,\real\lambda_2}{(x\,\|\vc v_1\|)^2+(y\,\|\vc v_2|)^2}, \\
|\pv(\vc w)|^2
=\frac{\|\varphi(\vc w)\|^2}{\|\vc w\|^2}
=\frac{(x^2\,\|\vc v_1\|)^2|\lambda_1|^2+(y\,\|\vc v_2\|)^2|\lambda_2|^2}{(x\,\|\vc v_1\|)^2+(y\,\|\vc v_2|)^2}. 
\end{gather*}
$\real\pv(\vc w)$ divides $\real\lambda_1$ and $\real\lambda_2$ into the ratio $(y\|\vc v_2\|)^2:(x\|\vc v_1\|)^2$. 
And $|\pv(\vc w)|^2$ divides $|\lambda_1|^2$ and $|\lambda_2|^2$ into the same ratio.
Since each point on $\ell(W)$ lies between $\lambda_1$ and $\lambda_2$, 
$\ell(W)$ is the geodesic with the endpoints $\lambda_1$ and $\lambda_2$. 
}

\thmenv{proposition}
{\label{loop-filled}
If $\dim V\ge3$, 
then $\D(W)\subset\Psi(\varphi)$ for any $2$-dimensional subspace $W\subset V$. 
}

\thmenv{proof}
{
It is enough to show the assertion in the case that $\D(W)$ has an interior point. 
In this case, $\ell(W)$ is a closed simple curve. 
Let $U$ be a $3$-dimensional subspace of $V$ including $W$. 
Put $\Proj^+(U)=U^\times/\R_{>0}$. 
$\Proj^+(U)$ is the quotient set of $U^\times$ identified by multiplying by positive real numbers, namely the set of half-lines starting at the origin. 
$\Proj^+(U)$ is homeomorphic to a spherical surface, and double covering of the projective plane $\Proj(U)$. 
The map $\pv$ is also regarded as a map on $\Proj^+(U)$. 
We have that $\pv(U^\times)=\pv(\Proj(U))=\pv(\Proj^+(U))$.  
For any $2$-dimensional subspace $W_1\subset U$, the quotient space $\Proj^+(W_1)=W_1^\times/\R_{>0}$ is a great circle on a spherical surface $\Proj^+(U)$. 

Let $\{\vc v_0,\,\vc v_1,\,\vc v_2\}$ be a basis of $U$ satisfying that $\{\vc v_0,\,\vc v_1\}$ is a basis of $W$. 
Put $\lambda_0=\pv(\vc v_0)$, $\lambda_1=\pv(\vc v_1)$, $\lambda_2=\pv(\vc v_2)\in\overline{\Hu}$, and $T=\langle\vc v_1,\vc v_2\rangle$. 
For any $\vc t\in\Proj^+(T)$, we define the closed oriented curve $c(\vc t)$ by 
\[c(\vc t):[0,\,1]\ni\theta\mapsto \pv(\cos(\pi\,\theta)\vc v_0+\sin(\pi\,\theta)\vc t))\in\ell(\langle\vc v_0,\vc t\rangle)\subset\pv(U^\times)\subset\Psi(\varphi).\]
$c(\vc t)$ is a closed curve passing through $\lambda_0$, 
because $c(\vc t)(0)=c(\vc t)(1)=\pv(\pm\vc v_0)=\lambda_0$. 
$\{c(\vc t)\}_{\vc t}$ is the continuous deformation of a closed curve on the complex plane. 
$c(\vc v_1)$ and $c(-\vc v_1)$ are the same curve $\ell(W)$, but have different orientations. 
Through the continuous deformation $\{c(\vc t)\}_{\vc t}$, the closed curve $\ell(W)$ is superimposed on itself in the opposite orientation. 
Therefore, the continuous deformation passes through all interior points of $\D(W)$. 
Any interior point of $\D(W)$ lies on a curve $c(\vc t)$ for some $\vc t\in\Proj^+(T)$. 
So, $\D(W)\subset\Psi(\varphi)$ holds. 
}

\thmenv{mainthm}
{
Leaf is convex on the Poincar\'e metric, for any endomorphism of dimension at least $3$ 
}

\thmenv{proof}
{
Let $\varphi$ be an endomorphism on a real inner product space $V$ of dimension at least $3$. 
Let $\lambda_1$, $\lambda_2\in\Psi(\varphi)$, $\vc v_1$, $\vc v_2\in V$ with $\lambda_1=\pv(\vc v_1)$ and $\lambda_2=\pv(\vc v_2)$, and $s(\lambda_1,\,\lambda_2)$ be  the geodesic connecting $\lambda_1$ to $\lambda_2$. 
By Proposition \ref{loop}, $s(\lambda_1,\,\lambda_2)$ is included in $\D(\langle\vc v_1,\,\vc v_2\rangle)$. 
By Proposition \ref{loop-filled}, $\D(\langle\vc v_1,\,\vc v_2\rangle)$ is included in $\Psi(\varphi)$. 
Hence, the leaf is convex on the Poincar\'e metric. 
}

The last theorem leads that the leaf of an endomorphism of dimension at least $3$ is simply connected, 
so the following theorem holds by Proposition \ref{leaf-simp-eigen}. 
It is interesting to note that the following algebraic property are derived from the geometric property. 

\thmenv{mainthm}
{
Leaf contains all eigenvalues whose imaginary part is nonnegative, for any endomorphism of dimension at least $3$ 
}

\section{Leaf and the eigenvalue geodesic polygon}

For any endomorphism on the real inner product space of finite dimension, the smallest convex domain in $\overline{\Hu}$ containing all eigenvalues with nonnegative imaginary part is a filled geodesic polygon whose vertexes are eigenvalues with nonnegative imaginary part. 
The domain is called the eigenvalue geodesic polygon. 
By the structure theorems (Theorems \ref{convex} and \ref{leaf-eigen}), the leaf includes the eigenvalue geodesic polygon, for each endomorphism on a real inner product space of dimension at least $3$. 
In generally, the leaf is larger than the eigenvalue geodesic polygon. 
We will show that they coincide for each normal endomorphism (Theorem \ref{normal-polygon}) in this section. 

Let $\varphi$ be an endomorphism on a real inner product space $V$ of finite dimension. 

\thmenv{proposition}
{\label{eigenSp}
Let $\lambda\in\Hu$ be an eigenvalues of $\varphi$ with positive imaginary part, and $\vc v\in V_\C$ an eigenvector of $\varphi$ corresponding to $\lambda$. 
Put $W=\langle\real\vc v,\,\imag\vc v\rangle\subset V$ and $\ell(W)=\pv(W^\times)$. 
\begin{enumerate}[ $($\rm a$)$]
\item 
$W$ is a $2$-dimensional $\varphi$-stable subspace of $V$. 
Hence $\ell(W)$ is the subleaf. 
\item 
Put $\lambda^+=\pv(\real \vc v)$, $\lambda^-=\pv(\imag \vc v)$ and $\beta=\vc v\cdot\overline{\vc v}/\|\vc v\|^2$. 
Then 
\[|\lambda^\pm|^2=\frac{|\lambda|^2\pm\real(\lambda^2\beta)}{1\pm\real\beta}, \quad 
\real\lambda^\pm=\frac{\real\lambda\pm\real(\lambda\beta)}{1\pm\real\beta}, \quad 
\imag\lambda^\pm=\imag\lambda\,\frac{\sqrt{1-|\beta|^2}}{1\pm\real\beta}.\]
\item  
$\lambda^\pm=\lambda$ holds, if and only if $\vc v\perp\overline{\vc v}$. 
\item 
If $\vc v\perp\overline{\vc v}$, then $\ell(W)=\{\lambda\}$. 
\item 
Three points $\lambda$, $\lambda^+$ and $\lambda^-\in\Psi(\varphi)$ lies on the same geodesic. 
\item 
$\ell(W)$ is a geodesic circle $\gC(\lambda;\,\frac{1-\sqrt{1-|\beta|^2}}{|\beta|})$. 
(rem. $\frac{1-\sqrt{1-|\beta|^2}}{|\beta|}\to0$ as $\beta\to0$) 
\end{enumerate}
}

\thmenv{proof}
{
\thmNum The assertion is already shown in \S 2. 
\\
\thmNum 
Since $\real\vc v=(\vc v+\overline{\vc v})/2$, $\imag\vc v=(\vc v-\overline{\vc v})/(2\unitCpx)$, and $\overline{\vc v}$ is an eigenvector of $\varphi$ corresponding to $\overline{\lambda}$, we have that 
\begin{alignat*}{2}
\|\real\vc v\|^2 
& = \frac{\|\vc v\|^2}{2}\,(1+\real\beta), & \qquad 
\|\imag\vc v\|^2 
& = \frac{\|\vc v\|^2}{2}\,(1-\real\beta), \\
\|\varphi(\real\vc v)\|^2 
& = \frac{\|\vc v\|^2}{2}\,(|\lambda|^2+\real\lambda^2\beta), & \qquad 
\|\varphi(\imag\vc v)\|^2 
& = \frac{\|\vc v\|^2}{2}\,(|\lambda|^2-\real\lambda^2\beta), \\
\real\vc v\cdot\varphi(\real\vc v) 
& = \frac{\|\vc v\|^2}{2}\,(\real\lambda+\real\lambda\beta), & \qquad 
\imag\vc v\cdot\varphi(\imag\vc v) 
& = \frac{\|\vc v\|^2}{2}\,(\real\lambda-\real\lambda\beta), 
\end{alignat*}
and then 
\begin{alignat*}{2}
|\lambda^+|^2 = \|\pv(\real\vc v)\|^2 
& = \frac{|\lambda|^2+\real\lambda^2\beta}{1+\real\beta}, & \qquad 
\real\lambda^+ = \real\pv(\real\vc v) 
& = \frac{\real\lambda+\real\lambda\beta}{1+\real\beta}, \\
|\lambda^-|^2 = \|\pv(\imag\vc v)\|^2 
& = \frac{|\lambda|^2-\real\lambda^2\beta}{1-\real\beta}, & \qquad 
\real\lambda^- = \real\pv(\imag\vc v) 
& = \frac{\real\lambda-\real\lambda\beta}{1-\real\beta}. 
\end{alignat*}
The imaginary parts of $\lambda^\pm$ ($>0$) are given by $\sqrt{|\lambda^\pm|^2-(\real\lambda^\pm)^2}$. 
\[(\imag\lambda^\pm)^2=|\lambda^\pm|^2-(\real\lambda^\pm)^2=\frac{|\lambda|^2\pm\real\lambda^2\beta}{1\pm\real\beta}-(\frac{\real\lambda\pm\real\lambda\beta}{1\pm\real\beta})^2=\frac{(\imag\lambda)^2\,(1-|\beta|^2)}{(1\pm\real\beta)^2}. \]
Since Cauchy-Schwarz inequality leads $|\beta|\le1$ and $|\real\beta|\le1$, we have that 
\[\imag\lambda^\pm=\imag\lambda\,\frac{\sqrt{1-|\beta|^2}}{1\pm\real\beta}. \]
\thmNum By the above formulae, 
$|\lambda^\pm|=|\lambda|$ if and only if $\lambda\beta\in\R$ and $\lambda\beta\le0$. 
$\real\lambda^\pm=\real\lambda$ if and only if $\beta\in\R$. 
Hence, $\lambda^\pm=\lambda$ if and only if $\beta=0$, i.e. $\vc v\perp\overline{\vc v}$. 
\\
\thmNum $\lambda$ and $\overline{\lambda}$ are the eigenvalues of the restriction $\varphi|_W$ of $\varphi$ on the two-dimensional $\varphi$-stable subspace $W$. 
By Propositions \ref{leaf-subsp-dim2} and \ref{leaf-dim2}, the closed curve $\ell(W)$ is a circle (or a single point set), and $\lambda$ lies on the circumference or the inside of $\ell(W)$. 
By \S 2, $\real\lambda\,\mtE+\imag\lambda\,\mtJ$ is the matrix representation of $\varphi|_W$ with respect to the basis $\real\vc v$, $\imag\vc v$ of $W$. 
Since $\vc v$ and $\overline{\vc v}$ are orthogonal, $\real\vc v$ and $\imag\vc v$ are also orthogonal and equal in length. 
$\ell(W)$ is the circle with radius $0$ by Proposition \ref{leaf-dim2}, so the single point set $\{\lambda\}$. 
\\
\thmNum 
Each geodesic of $\Hu$ is a part of line parallel to the imaginary axis or a semicircle whose center is on the real axis. 
The geodesic connecting through $\lambda^+$ and $\lambda^-$ is a part of line if $\real\lambda^+=\real\lambda^-$, or a part of semicircle if $\real\lambda^+\neq\real\lambda^-$. 

$\bullet$ In the case that $\real\lambda^+=\real\lambda^-$. 
It follows from (b) that 
\begin{gather*}
(\real\lambda^+=)\ \frac{\real\lambda+\real\lambda\beta}{1+\real\beta}=\frac{\real\lambda-\real\lambda\beta}{1-\real\beta}\ (=\real\lambda^-) \\
(\real\lambda+\real\lambda\beta)(1-\real\beta)=(\real\lambda-\real\lambda\beta)(1+\real\beta) \\
\therefore \, \real\lambda\,\real\beta=\real\lambda\beta
\end{gather*}
\[\therefore \, \real\lambda^\pm=\frac{\real\lambda\pm\real\lambda\beta}{1\pm\real\beta}=\frac{\real\lambda\pm\real\lambda\,\real\beta}{1\pm\real\beta}=\real\lambda. \]
$\lambda$ lies on the line which is the geodesic connecting through $\lambda^+$ and $\lambda^-$. 

$\bullet$ In the case that $\real\lambda^+\neq\real\lambda^-$. 
The geodesic connecting through $\lambda^+$ and $\lambda^-$ is a semicircle. 
Put $c\in\R$ be the center of the semicircle, and $d$ the radius. 
The calculation in the proof of Proposition \ref{leaf-dim2} leads that 
\[d^2-c^2=|\lambda^\pm|^2-2c\,\real\lambda^\pm \qquad c=\frac{|\lambda^+|^2-|\lambda^-|^2}{2\,(\real\lambda^+-\real\lambda^-)} \]
To show that $\lambda$ is on the semicircle, we need only show that $|\lambda-c|=d$ holds.
Since the equality $|\lambda-c|=d$ is equivalent to 
\[d^2=|\lambda-c|^2=|\lambda|^2-2c\,\real\lambda+c^2,\]
we show that $d^2-c^2=|\lambda|^2-2c\,\real\lambda$. 
\begin{align*}
& d^2-c^2 = |\lambda^\pm|^2-2c\,\lambda^\pm = |\lambda^\pm|^2-\frac{|\lambda^+|^2-|\lambda^-|^2}{\real\lambda^+-\real\lambda^-}\real\lambda^\pm \\
& \hspace{34pt}  = \frac{|\lambda^\pm|^2\,(\real\lambda^+-\real\lambda^-)-(|\lambda^+|^2-|\lambda^-|^2)\,\real\lambda^\pm}{\real\lambda^+-\real\lambda^-} 
 = \frac{|\lambda^-|^2\,\real\lambda^+-|\lambda^+|^2\,\real\lambda^-}{\real\lambda^+-\real\lambda^-} \\
& |\lambda|^2-2c\,\real\lambda =\frac{|\lambda|^2\,(\real\lambda^+-\real\lambda^-)-(|\lambda^+|^2-|\lambda^-|^2)\,\real\lambda}{\real\lambda^+-\real\lambda^-} 
\end{align*}
Substituting (b) into the right-hand side of the two equations above, both are equal to the following: 
\[\frac{2(|\lambda|^2\,\real\lambda\beta-\real\lambda\,\real\lambda^2\beta)}{(1+\real\beta)(1-\real\beta)(\real\lambda^+-\real\lambda^-)}.\]
Hence $|\lambda-c|=d$. 
Namely, $\lambda$ lies on the geodesic connecting through $\lambda^+$ and $\lambda^-$. 
\\
\thmNum 
In \S 3, we already see that the geodesic circle $\gC(\lambda;\,k)$ with center $\lambda=p+q\unitCpx$ for some $0\le k\le1$ is the circle with center $p+\frac{1+k^2}{1-k^2}q\unitCpx$ and radius  $\frac{2\,k}{1-k^2}q$. 
We give $k\in\R$ such that $\ell(W)=\gC(\lambda;\,k)$. 
In the case $\vc v\perp\overline{\vc v}$, we already get $\ell(W)=\{\lambda\}$ in (c). 
Since $\frac{1-\sqrt{1-\beta^2}}{|\beta|}\to0$ as $\beta\to0$, we consider that $\gC(\lambda;\,\frac{1-\sqrt{1-\beta^2}}{|\beta|})=\gC(\lambda;\,0)=\{\lambda\}$ for $\beta=0$. 
$\vc v\perp\overline{\vc v}$ means $\beta=0$, so the assertion holds in this case. 

Assume that $\vc v\not\perp\overline{\vc v}$, i.e. $\beta\neq0$. 
Let $\matrixform{cc}{a&b\\c&d}$ be a matrix representation of $\varphi|_W$ with respect to an orthonormal basis of $W$. 
Then $\lambda=\frac{(a+d)+\sqrt{D}}{2}$, where $D=(a+d)^2-4(ad-bc)$. 
By Proposition \ref{leaf-dim2}, the subleaf $\ell(W)=\Psi(\varphi;\,W)$ corresponding to $W$ is the circle with center $\frac{(a+d)+|b-c|\sqrt{-1}}{2}$ and radius $\frac{\sqrt{(a-d)^2+(b+c)^2}}{2}$. 
For the centers of circles $\gC(\lambda;\,k)$ and $\ell(W)$ to coincide, $p$, $q$ and $k$ satisfy that 
\[p=\frac{a+d}{2}, \qquad q=\frac{\sqrt{-D}}{2}, \qquad \frac{1+k^2}{1-k^2}\,q=\frac{|b-c|}{2}.\]
And then the radiuses of the circles coincide as the following. 
\[(\frac{2\,k}{1-k^2}\,q)^2=((\frac{1+k^2}{1-k^2})^2-1)\,q^2 = (\frac{|b-c|}{2})^2-(\frac{\sqrt{-D}}{2})^2
=(\frac{\sqrt{(a-d)^2+(b+c)^2}}{2})^2.\]
Hence $\ell(W)=\gC(\lambda;\,k)$. 
$\lambda^+$ lies on the geodesic circle $\gC(\lambda;\,k)$, so we have that 
\[k^2=\delta(\lambda,\,\lambda^+)^2=\frac{|\lambda-\lambda^+|^2}{|\lambda-\overline{\lambda^+}|^2}
=\frac{(\real\lambda-\real\lambda^+)^2+(\imag\lambda-\imag\lambda^+)^2}{(\real\lambda-\real\lambda^+)^2+(\imag\lambda+\imag\lambda^+)^2}
=\frac{1-\sqrt{1-|\beta|^2}}{1+\sqrt{1-|\beta|^2}}.
\]
Hence, 
\[\sqrt{1-|\beta|^2}=\frac{1-k^2}{1+k^2}\qquad 
\therefore\ |\beta|=\frac{2k}{1+k^2}. \]
Since $k$ is the root of the quadratic equation $|\beta|\,k^2-2\,k+|\beta|=0$ lying on the interval $[0,\,1]$, we obtain that 
\[k=\frac{1-\sqrt{1-|\beta|^2}}{|\beta|}.\] 
}

\thmenv{corollary}
{
\label{eigen2}
Let $\lambda$ be an eigenvalue of $\varphi$ with positive imaginary part. 
If there exist two $\C$-linearly independent eigenvectors corresponding to $\lambda$, then $\lambda$ belongs to $\Psi(\varphi)$. 
}

\thmenv{proof}
{
Let $\vc v_1$, $\vc v_2\in V_\C$ be $\C$-linearly independent eigenvectors corresponding to $\lambda$. 
If $\vc v_1\perp\overline{\vc v_1}$, then the last proposition (c) leads that $\pv(\real\vc v_1)=\lambda$, i.e. $\lambda\in\Psi(\varphi)$. 
Assume that $\vc v_1\not\perp\overline{\vc v_1}$. 
Put $\vc v=t\,\vc v_1+\vc v_2\in V_\lambda$ ($t\in\C$). 
The quadratic equation $\vc v\cdot\overline{\vc v}=(\vc v_1\cdot\overline{\vc v_1})\,t^2+2\,\real(\vc v_1\cdot\overline{\vc v_2})\,t+(\vc v_2\cdot\overline{\vc v_2})=0$ has a root $t_0\in\C$, because the coefficient of the leading term is nonzero. 
Put $\vc v_0=t_0\vc v_1+\vc v_2$, and then $\vc v_0\perp\overline{\vc v_0}$. 
The last proposition (c) leads that $\pv(\real\vc v_0)=\lambda$, i.e. $\lambda\in\Psi(\varphi)$. 
}

The last corollary is included in Theorem \ref{leaf-eigen}. 
By Theorem B, there exists a vector whose value of the function $\pv$ is the given eigenvalue. 
However, by the last corollary, one such vector is explicitly given by eigenvectors, if the eigenspace is of dimension at least $2$. 

Let $\lambda\in\overline{\Hu}$ be an eigenvalue of $\varphi$ with nonnegative imaginary part, and $V_\lambda\subset V_\C$ the eigenspace corresponding to $\lambda$. 
The real part $\real V_\lambda=\{\real\vc v\mid\vc v\in V_\lambda\}$ of $V_\lambda$ is a subspace of $V$. 
It is easy to show that $\real V_\lambda=(V_\lambda+\overline{V_\lambda})\cap V=\imag V_\lambda$ and $\overline{V_\lambda}=V_{\overline{\lambda}}$. 
Moreover, since $\real V_\lambda$ is $\varphi$-stable, 
$\pv((\real V_\lambda)^\times)$ is the subleaf $\Psi(\varphi;\,\real V_\lambda)$ corresponding to $\real V_\lambda$. 
In Proposition \ref{eigenSp}, we use the closed curve $\ell(\langle\real\vc v,\,\imag\vc v\rangle)$ for some $\vc v\in V_\lambda$, to see $\Psi(\varphi)$ around $\lambda$. 
We need the subleaf $\Psi(\varphi;\,\real V_\lambda)$ to look at the relationship between $\lambda$ and the other point on $\Psi(\varphi)$.  

\thmenv{proposition}
{\label{1-point}
Let $\lambda\in\overline{\Hu}$ be an eigenvalue with non-negative imaginary part. 
\begin{enumerate}[ $($\rm a$)$]
\item 
If $\lambda\in\R$, then 
$\Psi(\varphi;\,\real V_\lambda)=\{\lambda\}$. 
\item 
Assume that $\lambda\not\in\R$. 
$\Psi(\varphi;\,\real V_\lambda)=\{\lambda\}$ holds, if and only if $V_\lambda\perp \overline{V_\lambda}$. 
\item 
Assume that $\lambda\not\in\R$ and $V_\lambda\not\perp\overline{V_\lambda}$. 
Then 
$\Psi(\varphi;\,\real V_\lambda)$ is a geodesic circle with center $\lambda$ when $\dim \real V_\lambda=2$, and a geodesic disc (a filled geodesic circle) with center $\lambda$ when $\dim \real V_\lambda>2$. 
\end{enumerate}
}

\thmenv{proof}
{
\thmNum Since $\lambda\in\R$, $\real V_\lambda$ is the eigenspace in $V$ corresponding to $\lambda$. 
It follows from the proof of Proposition \ref{leaf-real-eigen} that $\Psi(\varphi;\,\real V_\lambda)=\pv((\real V_\lambda)^\times)=\{\pv(\vc v)\mid\vc v\in (\real V_\lambda)^\times\}=\{\lambda\}$. 
\\
\thmNum 
Assume that $V_\lambda\perp\overline{V_\lambda}$. 
Each $\vc v_1\in V_\lambda$ is orthogonal to its conjugate. 
Proposition \ref{eigenSp} (c) leads that $\pv(\real\vc v_1)=\lambda$. 
Then, $\Psi(\varphi;\,\real V_\lambda)=\pv((\real V_\lambda)^\times)=\{\lambda\}$. 

Conversely, we assume that $V_\lambda\not\perp\overline{V_\lambda}$. 
Take $\vc v_1$, $\vc v_2\in V_\lambda$ satisfying that $\vc v_1\not\perp\overline{\vc v_2}$. 
At least one of $\vc v_1$, $\vc v_2$ and $\vc v_1+\vc v_2$ is not orthogonal to the conjugate of itself. 
Because, if $\vc v_1\perp\overline{\vc v_1}$ and $\vc v_2\perp\overline{\vc v_2}$, then $(\vc v_1+\vc v_2)\not\perp\overline{\vc v_1+\vc v_2}$ as the following: 
\[(\vc v_1+\vc v_2)\cdot(\overline{\vc v_1+\vc v_2})
=\vc v_1\cdot\overline{\vc v_1}+\vc v_1\cdot\overline{\vc v_2}+\vc v_2\cdot\overline{\vc v_1}+\vc v_2\cdot\overline{\vc v_2}
=\vc v_1\cdot\overline{\vc v_2}+\vc v_2\cdot\overline{\vc v_1}
=2\,\vc v_1\cdot\overline{\vc v_2}\neq0.\]
Put $\vc v_0$ one of $\vc v_1$, $\vc v_2$, $\vc v_1+\vc v_2$ which is not orthogonal to the conjugate of itself. 
Proposition \ref{eigenSp} (c) leads that $\pv(\real\vc v_0)\neq\lambda$. 
Hence, $\Psi(\varphi;\,\real V_\lambda)=\pv((\real V_\lambda)^\times)\neq\{\lambda\}$. 
\\
\thmNum 
If $\dim \real V_\lambda=2$, then the subleaf $\Psi(\varphi;\,\real V_\lambda)$ is a geodesic circle with center $\lambda$. 
Assume that $\dim V_\lambda\ge3$. 
Then the subleaf $\Psi(\varphi;\,\real V_\lambda)$ is simply connected. 
By Proposition \ref{eigenSp} (f), $\ell(\langle\real\vc v,\,\imag\vc v\rangle)$ is a geodesic circle with center $\lambda$ for any $\vc v\in V_\lambda$. 
Since $\real V_\lambda=\bigcup_{\vc v\in V_\lambda}\langle\real \vc v,\,\imag\vc v\rangle$, we see that $\Psi(\varphi;\,\real V_\lambda)
=\bigcup_{\vc v\in V_\lambda}\ell(\langle\real\vc v,\,\imag\vc v\rangle)$, the union of geodesic circle with same center. 
The simply connected domain $\Psi(\varphi;\,\real V_\lambda)$ is a geodesic disc with center $\lambda$. 
}

\thmenv{remark}
{
The geometric property (Theorem \ref{convex}) was used in the above proof of the last proposition (c). 
However, it can be shown without it. 
The subleaf $\Psi(\varphi;\,\real V_\lambda)$ is the union of concentric geometric circles $\ell(\langle\real\vc v,\,\imag\vc v\rangle)$ ($\vc v\in\real V_\lambda$) with center $\lambda$, by Proposition \ref{eigenSp} (f). 
Moreover, Proposition \ref{eigenSp} (f) and the argument of the proof of Corollary \ref{eigen2} lead that the set of radiuses of the geometric circles $\ell(\langle\real\vc v,\,\imag\vc v\rangle)$ is a closed interval from $0$ to some nonnegative real number. 
So, the leaf $\Psi(\varphi;\,\real V_\lambda)$ is a geodesic disc with center $\lambda$. 
}

\thmenv{definition}
{
For any subset $A$, $B\subset\overline{\Hu}$, the union of all geodesics connecting a point on $A$ to a point on $B$ is called the {\it geodesic bridge} connecting $A$ to $B$, denoted by $A\asymp B$. 
}

The geodesic bridge $\asymp$ is a commutative semigroup on the power set of $\overline{\Hu}$, that is, it satisfies the associativity: $(A\asymp B)\asymp C=A\asymp(B\asymp C)$, and the commutativity: $A\asymp B=B\asymp A$. 
The geodesic bridge connecting finite convex subsets is the smallest convex set including all those subsets. 
Hence, the eigenvalue geodesic polygon of an endomorphism $\varphi$ is represented by $\{\lambda_1\}\asymp\cdots\asymp\{\lambda_r\}$ where $\lambda_1,\,\ldots,\,\lambda_r$ are the eigenvalue of $\varphi$ with nonnegative imaginary part. 

\thmenv{proposition}
{\label{leaf-subleaf}
Let $\varphi$ be an endomorphism on a real inner product space $V$, 
and $W_1$ and $W_2$ two nonzero $\varphi$-stable subspaces of $V$ satisfying that $W_1\perp W_2$. 
Then \[\Psi(\varphi;\,W_1\oplus W_2)=\Psi(\varphi;\,W_1)\asymp\Psi(\varphi;\,W_2).\] 
}

\thmenv{proof}
{
$\dim(W_1+W_2)\le2$ leads $\dim W_1=\dim W_2=1$, so $\Psi(\varphi;\,W_1)$ and $\Psi(\varphi;\,W_2)$ are single point sets. 
The elements $\lambda_1$ and $\lambda_2\in\R$ of these sets are eigenvalues of $\varphi|_{W_1}$ and $\varphi|_{W_2}$, respectively, 
Since $W_1\perp W_2$, the representation matrix of $\varphi|_{W_1\oplus W_2}$ of an orthonormal bases is the diagonal matrix whose components of diagonal are $\lambda_1$ and $\lambda_2$. 
By Proposition \ref{leaf-dim2}, the subleaf $\Psi(\varphi;\,W_1\oplus W_2)$ is a geodesic connecting through $\lambda_1$ and $\lambda_2$, that is the geodesic bridge connecting $\{\lambda_1\}$ to $\{\lambda_2\}$. 
Hence, $\Psi(\varphi;\,W_1\oplus W_2)=\Psi(\varphi;\,W_1)\asymp\Psi(\varphi;\,W_2)$. 

Assume that $\dim(W_1\oplus W_2)\ge3$. 
The subleaf $\Psi(\varphi;\,W_1\oplus W_2)$ is convex, so $\Psi(\varphi;\,W_1\oplus W_2)$ includes $\Psi(\varphi;\,W_1)\asymp\Psi(\varphi;\,W_2))$. 
We will show the inverse inclusion relation. 
By looking at the end points of the geodesics, $\Psi(\varphi;\,W_1)$ and $\Psi(\varphi;\,W_2)$ are included in $\Psi(\varphi;\,W_1\oplus W_2)$. 
Let $\lambda\in\Psi(\varphi\,W_1\oplus W_2)\smallsetminus(\Psi(\varphi;\,W_1)\cup\Psi(\varphi;\,W_2))$. 
Take $\vc v=\vc w_1+\vc w_2\in W_1\oplus W_2$ such that $\lambda=\pv(\vc v)$, $\vc w_1\in W_1$ and $\vc w_2\in W_2$. 
Since $\vc w_1$ and $\vc w_2$ satisfy the orthogonal condition in Proposition \ref{geod}, the curve $\ell(\langle\vc w_1,\,\vc w_2\rangle)$ is the geodesic connecting $\pv(\vc w_1)$ to $\pv(\vc w_2)$. 
Since $\lambda=\pv(\vc w_1+\vc w_2)\in\ell(\langle\vc w_1,\,\vc w_2\rangle)$, $\pv(\vc w_1)\in\Psi(\varphi;\,W_1)$ and $\pv(\vc w_2)\in\Psi(\varphi;\,W_2)$, $\lambda$ lies on the geodesic bridge connecting $\Psi(\varphi;\,W_1)$ to $\Psi(\varphi;\,W_2)$. 
Hence, 
\[\Psi(\varphi;\,W_1\oplus W_2)=\Psi(\varphi;\,W_1)\asymp\Psi(\varphi;\,W_2).\]
}

\thmenv{mainthm}
{
Leaf is the eigenvalue geodesic polygon, for any normal endomorphism of dimension at least 3, 
}

\thmenv{proof}
{
Let $\varphi$ be a normal endomorphism on a real inner product space $V$ of finite dimension. 
Then, eigenvectors of $\varphi$ corresponding to different eigenvalues are orthogonal. 
$V$ is the orthogonal sum of real parts of eigenspaces of all eigenvalues with nonnegative imaginary part. 
Let $\Lambda=\{\lambda_1,\,\ldots,\,\lambda_r\}$ be the set of all eigenvalues with nonnegative imaginary part, $\overline{\Lambda}$ the set of all eigenvalues, and $V_\lambda$ the eigenspace corresponding to $\lambda\in\overline{\Lambda}$. 
Since $\varphi$ is normal, $V_\lambda\perp V_{\lambda'}$ for any $\lambda$, $\lambda'\in\overline{\Lambda}$ with $\lambda\neq\lambda'$. 
And then $\Psi(\varphi;\,\real V_\lambda)=\{\lambda\}$ for any $\lambda\in\Lambda$ by Proposition \ref{1-point}, and $V$ is the $\varphi$-stable orthogonal sum of $\real V_{\lambda_1}$, $\cdots$, $\real V_{\lambda_r}$. 
Namely, $V=\real V_{\lambda_1}\oplus\cdots\oplus\real V_{\lambda_r}$, $\real V_{\lambda_i}\perp\real V_{\lambda_j}$ ($i\neq j$), and each $\real V_{\lambda_i}$ is $\varphi$-stable. 
By Proposition \ref{leaf-subleaf}, 
\[\Psi(\varphi)
=\Psi(\varphi;\,\real V_{\lambda_1})\asymp\cdots\asymp\Psi(\varphi;\,\real V_{\lambda_r})  
=\{\lambda_1\}\asymp\{\lambda_2\}\asymp\cdots\asymp\{\lambda_r\}\]
}

\thmenv{remark}
{\label{nonnormal-polygon}
In general, the converse of the last theorem does not hold. 
For example, 
take the endomorphism $\varphi=(-1)\oplus(0)\oplus(1)\oplus((3\,\mtJ+\mtK)/4)$ on $\R^5$, 
where $\mtJ$ and $\mtK$ are quadratic square matrixes defined in Example \ref{eg-leaf}. 
$(3\,\mtJ+\mtK)/4$ is not normal, so $\varphi$ is also not. 
The eigenvalues with nonnegative imaginary part are $-1$, $0$, $1$ and $\ssqrt{-1/2}$. 
The filled geodesic triangle $\Delta$ with vertexes $-1$, $0$ and $1$ is surrounded by three semicircles with center $1/2$ and radius $1/2$, with center $-1/2$ and radius $1/2$, and with center $0$ and radius $1$. 
The eigenvalue $\lambda=\ssqrt{-1/2}$ lies on the interior of $\Delta$, so $\Delta$ is also the eigenvalue geodesic polygon of $\varphi$. 
Put $W=\real V_\lambda$, then we have that $\varphi|_W=(3\,\mtJ+\mtK)/4$ and $\varphi|_{W^\bot}=(-1)\oplus(0)\oplus(1)$. 
The subleaf $\Psi(\varphi;\,W)$ is a circle with center $\alpha=3\unitCpx/4$ and radius $1/4$ by Proposition \ref{leaf-dim2}, and the subleaf $\Psi(\varphi;\,W^\bot)=\Delta$ by the last theorem. 
We see that the circle $\Psi(\varphi;\,W)$ is included in $\Delta$, by calculating the distance between the center $\alpha$ of the circle and the edges of $\Delta$. 
Therefore, the leaf $\Psi(\varphi)$ is the eigenvalue geodesic polygon. 
$\varphi$ and $\varphi|_{W^\bot}$ have the same leaf. 
By replacing or increasing the direct sum factor appropriately, it is possible to obtain infinitely many endomorphisms with the same leaf.
}

\section*{Acknowlegments}

The author would like to thank Hiroaki Nakamura, Yu Yasufuku, Michihisa Wakui and Kanji Namba for several crucial remarks and valuable information on relevant researches. 
Special thanks are due to Hiroaki Nakamura for his support and suggestions on various aspects of the author's continued research.

\end{document}